\documentclass[12pt]{article} %%%msart}
\usepackage{amsfonts}
\usepackage{amsmath}
\def\A{{\mathbb A}}
\usepackage{amsthm}
\def\Z{{\mathbb Z}}
\def\C{{\mathbb C}}

\def\N{{\mathbb N}}
\def\CA{{\mathcal A}}

\DeclareMathSymbol\crossrt{\mathrel}{AMSb}{"6E}
\DeclareMathSymbol\crosslt{\mathrel}{AMSb}{"6F}

\def\ts{\otimes}

\newtheorem{lemma}{Lemma}[section]

\newtheorem{proposition}[lemma]{Proposition}
\newtheorem{corollary}[lemma]{Corollary}
\theoremstyle{definition}

\newtheorem{definition}[lemma]{Definition}

\newtheorem{remark}[lemma]{Remark}

\usepackage{times}
\begin{document}
\title{Twisted Hochschild Homology \\ of Quantum Hyperplanes}
\author{Andrzej Sitarz
\footnote{Partially supported by Polish State Committee for Scientific Research
(KBN) under grant 2 P03B 022 25} \\
%\address
{\em Institute of Physics, Jagiellonian University,}\\
{\em Reymonta 4, 30-059 Krak\'ow, Poland} \\
e-mail: sitarz@if.uj.edu.pl }
%%%%%%%%%%%%%%%%%%%%%%%%%%%%%%%%%%%%%%%%%%%%%%%%%%%%%%
%%%%%%%%%%%%%%%%%%%%%%%%%%%%%%%%%%%%%%%%%%%%%%%%%%%%%%
\maketitle
%%%%%%%%%%%%%%%%%%%%%%%%%%%%%%%%%%%%%%%%%%%%%%%%%%%%%%
\begin{abstract}
We calculate the Hochschild dimension of quantum
hyperplanes using the twisted Hochschild homology.
\end{abstract}
\section{Introduction}

In noncommutative geometry the Hochschild and cyclic homology take
the role of differential de Rham complexes and de Rham homology.
It was observed quite early that in many canonical noncommutative
examples the Hochschild dimension of a noncommutative space
differs from the one of its commutative limit. This seems to be
the feature of quantum deformations, in particular the quantum
$SU_q(2)$ group and its homogeneous spaces (Podles spheres), which
show the "dimension drop" \cite{Jap-all}. Similarly, in the case
of generic quantum deformations of hyperplanes the effect is even
bigger, as irrespectively of the rank their Hochschild dimension
is $1$ \cite{Wa93}.

The twisted (modular) Hochschild and cyclic homology was
introduced by Kustermans, Murphy and Tuset \cite{KMT} based on the
notion of modular automorphisms and its relation with the Haar
state on the algebra on the quantum group. The latter satisfies
{\em twisted} cyclicity, which could be extended to closed graded
twisted traces on the differential algebras over the quantum
group, the latter giving explicit presentation of twisted cyclic
cocyles. The notions of the twisted (modular) Hochschild and
cyclic homologies became more attractive after explicit
demonstration of twisted cyclic cocycles of dimension $3$ for
$SU_q(2)$ \cite{KMT} and of dimension $2$ for the standard Podles
sphere, as well as the demonstration of the Connes-Moscovici local
formula or this cocycle \cite{Schm03} using the spectral triple
derived in \cite{DoSi02}. More recent calculations of examples
\cite{Had,Kra-Had} of twisted Hochschild cohomology for quantum
spaces confirm that for certain automorphisms the dimension drop
does not occur.

Apart from the motivation by the quantum group symmetries and the
Haar state, the twisted Hochschild and cyclic complexes can be
introduced for any algebra automorphism. In this paper we shall
investigate the example of quantum hyperplanes.

\section{Preliminaries}

Throughout the paper we work over the field of complex numbers.
Let $\sigma$ be an automorphism of the algebra $\CA$ and let
$\CA_\sigma$ be the vector space $\CA$ as a $\CA \otimes \CA^{op}$
right module with the following right-$\CA$-module structure:
\begin{equation}
b \cdot (a_0, a_1) := \sigma(a_1) b a_0,
\label{sigbi}
\end{equation}

We can consider now three types of homologies (we use
three different names in order to distinguish between
them later).

{\bf Natural twisted Hochschild homology.} This, we define
as the usual Hochschild homology of $\CA$ with values in
the bimodule $\CA_\sigma$, $HH(\CA,\CA_\sigma)$. It is
the homology of the following Hochschild chain complex:
\begin{equation}
\cdots \xrightarrow{b}  C_{n+1} \xrightarrow{b}
C_{n} \xrightarrow{b} \cdots \xrightarrow{b} C_1
\xrightarrow{b} C_0,
\label{cc1}
\end{equation}
where $C_n = \CA_\sigma \ts \CA^{\ts (n)}$.

{\bf Invariant twisted Hochschild homology}. We can view the
automorphism $\sigma$ as the generator of the
action of $\Z$ on the algebra $\CA$. One may easily verify
that the map $b$ of the chain complex (\ref{cc1}) is
equivariant with respect to this action, that is, it commutes
with the action of the group. Therefore, if we take a subspace
containing only invariant elements
$C_n^{inv} = (\CA_\sigma \ts \CA^{\ts n})_{inv}$,
that is: \mbox{$ a_0 \ts a_1 \ts \cdots \ts a_n
\in C_n^{inv}$} if and only if:
$$ a_0 \ts a_1 \ts \cdots \ts a_n =
\sigma(a_0) \ts \sigma(a_1) \ts \cdots \ts \sigma(a_n). $$
then the following is a subcomplex of (\ref{cc1}):
\begin{equation}
\cdots \xrightarrow{b}  C_{n+1}^{inv} \xrightarrow{b}
C_{n}^{inv} \xrightarrow{b} \cdots \xrightarrow{b} C_1^{inv}
\xrightarrow{b} C_0^{inv},
\label{cc2}
\end{equation}
and we shall denote its homology by $H(\CA,\CA_\sigma)_{inv}$.

{\bf Twisted Hochschild homology.} Finally, we can consider the
quotient of the Hochschild complex by the image of the $1-\sigma$
map, $C^{n}/(1-\sigma)$. We obtain in this way another chain complex:
\begin{equation}
\cdots \xrightarrow{b}  C_{n+1}/(1-\sigma) \xrightarrow{b}
C_{n}/(1-\sigma) \xrightarrow{b} \cdots \xrightarrow{b} C_1/(1-\sigma)
\xrightarrow{b} C_0/(1-\sigma),
\label{cc3}
\end{equation}

and we shall denote its homology by $HH_\sigma(\CA)$ and call it
(after \cite{KMT}) the twisted Hochschild homology of $\CA$. Since
only for such chains twisted cyclicity condition makes sense, this
complex is the one related to twisted cyclic homology.

Later, we shall need the following lemma.

\begin{lemma}\label{eqlemma}
Let $(K_*,d)$ and $(C_*,b)$ be chain complexes with the
defined action of a group $G$ such that both $b$ and $d$
are equivariant maps. If $f: K_* \to C_*$ is an equivariant
quasi-isomorphism of chain complexes and:
\begin{equation}
b^{-1}(K^{inv}_n \cap B(K_n)) = K^{inv}_{n+1} + Z(K_{n+1}),
\label{cond1}
\end{equation}
then the restriction of $f$ to the invariant subcomplexes
is also a quasi-isomorphisms.
\end{lemma}

\begin{proof}
Clearly $f: K_*^{inv} \to C_*^{inv}$ is a morphism of
chain complexes, so we need to proof that the induced map
in their homology is an isomorphism. Take \mbox{$k \in K_n^{inv}$}
such that $bk=0$ and assume that the class in $H(C_*^{inv})$
of its image $f(k)$ is trivial. Clearly, then the class
of $f(k)$ in $H(C_*)$ is trivial and since $f_*$ is
a quasi-isomorphism we know that the class of $k$
in $H(K_*)$ is trivial, that is that $k = b k_0$ for
some $k_0 \in K_{n+1}$. Now, using the assumption (\ref{cond1})
we know that $k_0$ could be presented as $k_0' + k'$, where
$k_0' \in K^{inv}_{n+1}$ and $k' \in Z(K_{n+1})$. Then
$b k_0' = k$ within $K_*^{inv}$ and the class of $k$ is
indeed trivial in $H(K_*^{inv})$.
\end{proof}

We shall end the section by introducing the definition
of a {\em scaling automorphism} and proving two important
lemmas on twisted Hochschild homology in that case.

\begin{definition}\label{scadef}
We say that $\sigma$ is a {\em scaling automorphism} of
the algebra $\CA$, if there exists a basis of $\CA$ as
a vector space $\{a_i\}_{i \in I}$ such that for each $i \in I$:
$\sigma(a_i) = p_i a_i$, $p_i \in \C, p_i \not=0$.
\end{definition}

Note that this property of $\CA$ extends naturally on
tensor products of $\CA$ with itself, so it remains true
for spaces of chains $C_n$, $n=0,1,2,\ldots$.

\begin{lemma}
For any scaling automorphism $\sigma$ (\ref{scadef}) the
invariant twisted Hochschild complex $C_*^{inv}$ is
isomorphic to the quotient complex $C_*/(1-\sigma)$.
\end{lemma}
\begin{proof}
It is sufficient to show that for each $n$,
$C_n = (1-\sigma) C_n \oplus C_n^{inv}$. Since the
boundary $b$ commutes with $\sigma$, then the natural
isomorphism, $C_n^{inv} \to C_n/(1-\sigma)$, which we have
for each $n=0,1,2,\ldots$, commutes with $b$.

Since, as we observed earlier, there exists a basis
$\{c_i\}$ of $C_n$, such that \mbox{$\sigma(c_i) = \gamma_i c_i$}
we easily see that $c_i$ is in the image of $(1-\sigma)$
unless $\gamma_i=1$ when it is $\sigma$-invariant.
\end{proof}

\section{The multidimensional quantum hyperplanes}

Let $V$ be a vector space of dimension $N$ and $T(V)$
its tensor algebra. If $q_{ij}$ is a complex matrix we denote
by $S_Q(V)$ the quantum symmetric algebra:
$$ S_Q(V) = T(V) / I_Q, $$
where $I_Q$ is the ideal in $T(V)$ generated by all elements
$$ X_{ij} = x_i \ts x_j - q_{ij} x_j \ts x_i, \;\;\; i < j,
 x_i,x_j \in V.$$

We always assume $q_{ii}=1$ for all $i=1,2,\ldots,N$ and
$q_{ji} =  q_{ij}^{-1}$ for $i < j$.

We introduce also the quantum antisymmetric algebra $\Lambda_Q(V)$
defined as $T(V)/ A_Q$, where $A_Q$ is the ideal in $T(V)$
generated by elements:
$$ A_{ij} = x_i \ts x_j + (q_{ij})^{-1} x_j \ts x_i, \;\;\; i < j,
 x_i,x_j \in V. $$
(Note that in case of an arbitrary field one needs some extra
generators to define correctly the antisymmetric algebra:
\mbox{$ A_{ii} = (1+ \hat{q}) x_i \ts x_i$}, where
$\hat{q} = q_{ij} q_{ji}$.)

\begin{remark}\label{scadef1}
Each scaling automorphism $\sigma$ of the quantum hyperplane
algebra is defined by its action on the generators,
for each $i=1,\ldots,N$:
$\sigma(x_i) = p_i x_i$, $p_i \in \C, p_i \not=0$.
By $S_Q^\sigma(V)$ we shall denote $S_Q(V)$ with the bimodule
structure (\ref{sigbi}).
\end{remark}

Following Wambst, we define:

\begin{definition}
The Koszul resolution of $S_Q(V)$ is the
following complex:
$$ K_n = S_Q(V) \ts \Lambda_Q^n(V) \ts S_Q(V), $$
with the differential $ d_n: K_n \mapsto K_{n-1}$:
$$
\begin{aligned}
&d_n ( a \ts x_{i_1} \wedge x_{i_2} \wedge \cdots
\wedge x_{i_n} \ts b ) = \\
&\phantom{xxx} = \sum_{k=1}^n (-1)^{k+1} \left(
(\prod_{s=1}^{k-1} q_{i_si_k}) a x_{i_k} \ts
x_{i_1} \wedge x_{i_2} \wedge \cdots \widehat{x_{i_k}} \cdots
\wedge x_{i_n} \ts b  \right. \\
&\phantom{xxx} \left. + (\prod_{s=k+1}^{n} q_{i_ki_s}) a  \ts
x_{i_1} \wedge x_{i_2} \wedge \cdots \widehat{x_{i_k}} \cdots
\wedge x_{i_n} \ts x_{i_k} b \right).
\end{aligned}
$$
In the formulae above $i_1 <  i_2 < \cdots < i_n$ and
$\widehat{x_{i_k}}$ denotes that this element falls out from the
product in $\Lambda^*_Q(V)$.

We denote this complex  \mbox{$K(S_Q(V),S_Q(V) \ts S_Q(V)^{op})$}.
\end{definition}

In his work Wambst proved (\cite{Wa93}, Proposition 4.1)
that this is a projective resolution over identity
of $S_Q(V)$ as a module over $S_Q(V) \ts S_Q(V)^{op}$.
Moreover, he constructed explicitly an $S_Q(V) - S_Q(V)$
linear morphism $\gamma'$ between this Koszul resolution
and the standard Hochschild resolution $C_*(S_Q(V))$,
thus enabling the direct calculation
of the Hochschild homology groups. We skip here detailed
presentation of the morphism referring to the mentioned
work (see \cite{Wa93}, Lemma 5.3).

Similarly as he did in the theorem 5.4, by applying
$S_Q^\sigma(V) \ts_{S_Q(V) \ts S_Q(V)^{op}}$ to both
resolutions, we get a quasi-isomorphism between the
"twisted" versions of Koszul and Hoschschild complexes.
Therefore the task of calculating the Hochschild
homology groups $HH(S_Q(V), S_Q^\sigma(V))$ is reduced
to calculation of the homology of the Koszul complex
$K_*(S_Q(V), S_Q^\sigma(V))$.

%%%%%%%%%%%%%%%%%%%%%%%%%%%%%%%%%%%%%%%%%%%%%%%%%%%%%%
Before we proceed with the calculation of respective
homology groups let us notice that the Hochschild
chain complex and the Koszul chain complex for the
quantum hyperplane and the scaling isomorphism
$\sigma$ satisfy the assumptions of the lemma (\ref{eqlemma}):

\begin{lemma}
For the Koszul chain complex of the quantum hyperplane
the condition (\ref{cond1}) is valid for the scaling
automorphism $\sigma$.
\end{lemma}
\begin{proof}
We shall prove a slightly more general statement,
keeping in mind that the Koszul chain complex has also
a basis, which is scaled by the automorphism $\sigma$.

Let $\{k_i\}_{i \in I}$ be such a basis of $K_*$,
that the automorphism $\sigma$ is scaling,
$\sigma(k_i) = \lambda_i k_i$, $\lambda_i \in \C$,
for each $i \in I$. Then for every $k \in K_{n+1}$
we have:
$$ k = \sum_i \alpha_i k_i,$$
where only finitely many of $\alpha_i$ are different from
$0$, so only finitely many possible scaling coefficients
$\lambda_i$ appear. Name the set of this coefficients
$\Lambda$ and set $k_\lambda$ the part of $k$, which has
the scaling coefficient $\lambda$. Then:
\begin{equation}
k = \sum_{\lambda \in \Lambda} k_{\lambda}. \label{eq5}
\end{equation}
Next, assume that $bk \in K_n^{inv}$. Since $b$ commutes
with the action of $\sigma$, the action of $\sigma$ on
each component $ b k_\lambda$ shall have the same scaling
coefficient $\lambda$. However, since $bk \in K_n^{inv}$, its
scaling coefficient must be $1$. On the other hand
$bk$ has parts $b k_\lambda$, with scaling coefficients
$\lambda \not=1$, so if their sum vanishes then all such
elements must be zero. Thus we obtain that in the sum
(\ref{eq5}) $k_1$ might be arbitrary but for all
$\Lambda \ni \lambda \not=1$ we must have $b k_\lambda=0$.
Therefore if $b k \in K_n^{inv}$ then
$k = k_1 + \sum_{\lambda \not= 1} k_\lambda$,
$k_1 \in K_{n+1}^{inv}$ and $k_\lambda \in Z(K_{n+1})$.
\end{proof}

{}From the explicit construction of the quasi-isomorphism
$\gamma$ in \cite{Wa93} we directly conclude:

\begin{remark}
The quasi-isomorphism $\gamma$ is equivariant with respect
to the action of $\Z$ set by the scaling automorphism
$\sigma$,
\end{remark}

It is the obvious conclusion from the construction of $\gamma$,
which involves only operations of (deformed) antisymmetrisations.

Finally we obtain:

\begin{corollary}
The homology groups of the Hochschild complex of $S_Q(V)$
with values in $S(Q)_\sigma$ are isomorphic with the
corresponding homology groups of the Koszul complex.
Similarly, invariant twisted Hochschild homology is isomorphic
to the homology of the invariant subcomplex of the Koszul
complex.
\end{corollary}

\pagebreak[3]

\subsection{Homology groups}

We begin with natural twisted Hochschild homology
of the quantum hyperplanes $HH(S_Q(V), S_Q(V)_\sigma)$.

We use the notation of \cite{Wa93}: let $\alpha \in \N^N$ denote
a natural multiindex and \mbox{$\beta \in \{0,1\}^N$}
a multiindex with values $0$ and $1$. We use the following
shorthand notation:
$$
\begin{aligned}
&x^\alpha = x_1^{\alpha_1} x_2^{\alpha_2} \cdots x_N^{\alpha_N},\\
&x^\beta = x_1^{\beta_1} \wedge x_2^{\beta_2} \wedge \cdots \wedge x_N^{\beta_N}.
\end{aligned}
$$

The length of $\alpha$ and $\beta$ defined as sum of all subindices is
denoted by $|\alpha|$, $|\beta|$, respectively. We introduce also
$\kappa_i \in \{0,1\}^N$, such that $\kappa_i(j)$ is zero for all $j \not= i$
and $\kappa_i(i)=1$.

Then, in the complex $S_Q(V)_\sigma \otimes_{S_Q(V) \ts S_Q(V)^{op}}
K(S_Q(V),S_Q(V) \ts S_Q(V)^{op})$ the differential on homogeneous
elements of the type $x^\alpha \ts x^\beta$ is:

$$
d ( x^\alpha \ts x^\beta ) =
\sum_{i=1}^N
\delta_1^{\beta(i)} \Omega(\alpha,\beta,i)
x^{\alpha+\kappa_i} \ts x^{\beta-\kappa_i}.
$$

The first factor $\delta_1^{\beta(i)}$ assures that the sum is only
over such indices $i$ for which $\beta(i)$ does not vanish and:

$$
\begin{aligned}
& \Omega(\alpha,\beta,i) = \\
& \phantom{xxx}
= (-1)^{ \stackrel{\stackrel{i-1}{\sum}}{s=1} \raise4pt\hbox{$\beta(s)$}} \left(
(\prod_{s=1}^{i-1} q_{si}^{\beta(s)}) (\prod_{r=i+1}^{N}
q_{ir}^{-\alpha(r)})
-p_i (\prod_{s=i+1}^{N} q_{is}^{\beta(s)})
 (\prod_{r=1}^{i-1} q_{ri}^{-\alpha(r)})
\right)
\end{aligned}
$$

By direct calculation it could be shown that $\Omega(\alpha,\beta,i)=0$
if and only if $x^{\alpha+\beta} x^i = \sigma(x^i) x^{\alpha+\beta}$.

Let us denote by $C_\sigma$ the set of all multiindices
$\gamma \in \N^N$, such that either $\gamma(i)=0$ or
$x^\gamma x^i = \sigma(x^i) x^\gamma$. From the earlier remark
we see that if $\alpha+\beta \in C_\sigma$ then
$d (x^\alpha \ts x^\beta )=0$. Since the differential $d$ preserves
the total multi-grade, $\alpha+\beta$, of the element
$x^\alpha \ts x^\beta$, to calculate the homology of the complex we
might restrict ourselves to the elements of a fixed multi-grade.

We have:
\begin{proposition}[compare theorem 6.1 \cite{Wa93}]\label{math}
The Hochschild homology of the quantum multiparameter
hyperplane $S_Q(V)$, with values in $S_Q^\sigma(V)$ is:
$$ HH_n(S_Q(V),S_Q^\sigma(V)) = \bigoplus_{\stackrel{\beta \in \{0,1\}^N}{|\beta|=n}}
\bigoplus_{\stackrel{\alpha \in \N^N}{\alpha + \beta \in C_\sigma}}
\C x^\alpha \ts x^\beta. $$
\end{proposition}

%%%%%%%%%%%%%%%%%%%%%%§§§§§§§§§§§§§§§§§§§§§§§§§§§§§§§§§§§§§§§

Proof follows exactly the proof of theorem 6.1 in \cite{Wa93}.
First, it is clear that all homogeneous elements $x^\alpha
\ts x^\beta$ for $\alpha+\beta \in C_\sigma$ are in the homology
groups of the complex. It remains to show that for
$\alpha+\beta \notin C_\sigma$ we might construct homotopy $h$
showing that this part of the complex is acyclic.

We define:
$$ h( x^\alpha \ts x^\beta ) =
\frac{1}{||\alpha+\beta||}
\sum_{i=1}^N \omega(\alpha,\beta,i)
x^{\alpha-\kappa_i} \ts x^{\beta+\kappa_i},
$$
where:
$$
\omega(\alpha,\beta,i) =
\begin{cases}
0 & \hbox{if\ }\; \alpha+\beta \in C_\sigma, \\
0 & \hbox{if\ }\; \beta(i)=1, \\
0 & \hbox{if\ }\; \alpha(i)=0,\\
\Omega(\alpha-\kappa_i,\beta+\kappa_i,i)^{-1} & \hbox{otherwise}.
\end{cases}
$$

Then one can calculate that:
$$ d h+h d = \hbox{id}.$$
For details of the proof we refer again to \cite{Wa93}.

Let us recall here the notion of {\em generic} deformation
parameters $q_{ij}$. We say that the parameters are not
generic if there exists such $\gamma \in \N^N$, $\gamma\not=0$
which is not $\kappa_j$ for some $j=1,\ldots,N$ and for
every $i$ such that $\gamma(i)>0$ we have:
$$ \prod_{k=1}^{N} (q_{ki})^{\gamma(k)} =1. $$

Note that this is equivalent to stating that the
center of any subalgebra generated within $S_Q(V)$
by more than two of its generators is nontrivial.

We are ready now to prove the main theorem:

\begin{proposition}
For each multiparameter quantum hyperplane $S_Q(V)$ with
generic deformation parameters $q_{ij}$ there
exist a scaling automorphism $\sigma$ such that the top
natural twisted Hochschild homology group of $S_Q(V)$ is
of dimension $N$, in particular, there exists a unique choice
of these parameters such that the top class in the twisted
Hochschild homology is given by $ 1 \ts x^1 \wedge x^2
\wedge \cdots \wedge x^N$ in the Koszul complex.
\end{proposition}

\begin{proof}
Clearly $HH_n(S_Q(V),S_Q^\sigma(V)) = 0$ for $n>N$. From the previous
proposition we need to find such automorphism
$\sigma$ that the space of elements $x^\alpha \ts x^\beta$ with
$|\beta|=N, \alpha+\beta \in C_\sigma$ were not empty. Since $\beta(i)=1$ for
all $i=1,2, \ldots,N$ (we denote it by $1_N$) we see that the
condition becomes:
$$ p_i x_i x^{\alpha+1_N} = x^{\alpha+1_N} x_i, \;\;\;
\forall i=1,2,\ldots,N.$$

Using the commutation rules for the $S_Q(V)$ we obtain:

$$ p_i \left( \prod_{j=1}^{i-1} (q_{ji})^{-\alpha(j)-1} \right)
=  \left( \prod_{j=i+1}^{N} (q_{ij})^{-\alpha(j)-1} \right).$$

and we might rewrite it as:

$$
p_i = \left( \prod_{j=1}^{i-1} (q_{ji})^{\alpha(j)+1} \right)
 \left( \prod_{j=i+1}^{N} (q_{ij})^{-\alpha(j)-1} \right)
= \left( \prod_{j=1}^{N} (q_{ji})^{\alpha(j)+1} \right).
$$

Therefore for any choice of $\alpha$ we might set the numbers
$p_i$ so that the highest twisted Hochschild homology group
does not vanish. In particular for the canonical choice
$\alpha=0$ we have:

\begin{equation}
p_i  = \left( \prod_{j=1}^{N} (q_{ji}) \right). \label{catw}
\end{equation}

Then $HH_N(S_Q(V),S_Q^\sigma(V)) = \C$ unless the family
of parameters $q_{ij}$ is not generic. We shall call such
scaling automorphisms as set by (\ref{catw}) {\em canonical
automorphisms}.
\end{proof}

Finally, we obtain:

\begin{proposition}\label{math2}
The twisted Hochschild homology of the quantum
hyperplane for the canonical scaling automorphism
(\ref{catw}) is, in dimension $n$:
$$ HH_n^\sigma(S_Q(V)) =
\bigoplus_{\stackrel{\beta \in \{0,1\}^N}{|\beta|=n}}
\bigoplus_{\stackrel{\alpha \in \N^N}{\alpha + \beta \in C_\sigma}}
\C x^\alpha \ts x^\beta. $$

where $C_\sigma$ is the set of all multiindices $\gamma \in \N^N$
such that for each $i=1,2,\ldots N$ either $\gamma(i)=0$ or
\begin{equation}
\prod_j q_{ji}^{\gamma(j)-1} = 1. \label{eqcond5}
\end{equation}
\end{proposition}

\begin{proof}
The proof is a direct consequence of the twisted Hochschild
homology construction, Proposition (\ref{math}) and the
formula (\ref{catw}) for the coefficients of the canonical
scaling automorphism.

First, notice that for the canonical scaling automorphism
the condition that $\gamma \in C_\sigma$ could be rewritten
(for $\gamma(i)\not=0$) exactly as in (\ref{eqcond5}).

It remains only to show that the elements of the Koszul
complex of the form $x^\alpha \ts x^\beta$, with
$\alpha+\beta \in C_\sigma$ are invariant with respect to
the canonical automorphism.

Let us verify explicitely the action of the canonical
scaling automorphism $\sigma$ on $x^{\gamma}$ for
$\gamma \in C_\sigma$:

$$
\sigma(x^{\gamma}) = \prod_i (p_i)^{\gamma(i)}
 = \prod_i \prod_j (q_{ji})^{\left(\gamma(i)\right)} =
\prod_i \prod_j (q_{ji})^{\left(\gamma(i)\gamma(j)\right)} = 1.
$$

where we have used first formula (\ref{catw}),
then (\ref{eqcond5}) and finally $q_{ij} q_{ji} = 1$.

Since we have shown already that for scaling automorphisms
the twisted Hochschild homology is equal to the invariant
twisted Hochschild homology we obtain the desired result.
\end{proof}

\subsection{The two-dimensional quantum plane}

In the simplest possible case of the generic two-dimensional
quantum plane $\A_q^2$ given by the relation $xy=qyx$, ($q$ not
a root of unity) we have the canonical scaling automorphism:
$$ \sigma(x) = \frac{1}{q} x, \;\;\;\;\; \sigma(y) =q y. $$

The twisted Hochschild homology groups of the quantum plane are:
$$
\begin{aligned}
& HH_0^\sigma(\A^2_q) = \C \oplus \C xy \\
& HH_1^\sigma(\A^2_q) =  (\C x \ts y) \oplus (\C y \ts x), \\
& HH_2^\sigma(\A^2_q) = \C 1 \ts x \wedge y, \\
& HH_n^\sigma(\A^2_q) = 0, \;\; n>3.
\end{aligned}
$$

where we have written the representatives from each of the
homology class from the (invariant) Koszul chain complex.

To show it, we use the Proposition (\ref{math2}). In dimension
$0$ we have $\beta \equiv 0$ and two possibilities for
$\alpha$, which belong to $C_\sigma^{inv}$: either
$\alpha(i)=0, i=1,2;$ or $\alpha(i)=1,i=1,2$. In dimension
$1$ we proceed likewise.

Finally, consider the example of the multidimensional
generalization of the quantum plane, the one-parameter
quantum hyperplane, which we use to illustrate that the
explicit calculations get quite technical in concrete examples.
The algebra, $\A_q^N$, is defined as an algebra generated by
$N$ generators $x_i$, $i=1,\ldots,N$ and relations (we
again assume that $q$ is not a root of unity):
\begin{equation}
x_i x_j  = q x_j x_i, \;\;\;
N \geq i > j \geq 1,
\end{equation}
The canonical scaling transformations (\ref{catw}) for $\A_q^N$ are:
\begin{equation}
\sigma(x_i) =  q^{N-2i+1} x_i, \;\;\; i=1,\ldots,N.
\end{equation}

The space $C_\sigma$ consists of all natural solutions to
the set of linear equations for $\gamma$:
\begin{equation}
P_k \left(\begin{array}{ccccc}
0 & 1 & 1 & \cdots & 1 \\
-1 & 0 & 1 &  \cdots & 1 \\
\vdots & & & & \vdots \\
-1 & \cdots & -1 & 0 & 1 \\
-1 & \cdots & -1 & -1 & 0
\end{array} \right) \gamma =
P_k \left(
\begin{array}{c}
N - 2 + 1 \\
N - 4 + 1 \\
\vdots \\
2 (N-1) - N + 1 \\
2 N - N + 1 \\
\end{array} \right),\label{eq1}
\end{equation}
where $P_k$ is any diagonal matrix with $0,1$ on the diagonal
(the matrix, which projects the equation on the non-zero values of
$\gamma$.

Using this technical description we have the tools to calculate the
twisted Hochschild homology in any dimension. In particular, in the
top dimension $n=N$, we have the only solution $\gamma = (1,1,\ldots,1)$.

\section{Conclusions}

We have demonstrated that for quantum hyperplanes there
exist automorphisms such that the top non-vanishing
twisted Hochschild homology group is of "classical"
dimension. This merely confirms that the twisted version
of Hochschild homology seems to be more adapted to the case
of quantum deformations. Note that in the discussed case we
did not start with a given automorphism but rather found
them (among some special type of scaling automorphisms)
using the requirement for the non-vanishing of the top
dimension.

Our results could be, of course, used also for the calculations
of the twisted cyclic homology. In fact, from the Connes long
exact sequence linking cyclic and Hochschild homology, which
is applicable also to the twisted case, we might immediately
get some results. In particular for $N$-dimensional quantum
hyperplanes the twisted cyclic homology stabilizes starting
from the dimension $N-1$. More detailed analysis of this
issue as well as of the case of non-generic deformation
parameters shall be treated elsewhere.

{\bf Acknowledgements:} The author thanks Ulrich Kr\"ahmer for 
discussions and comments on the manuscript.

\end{document}